# PERCOLATION ON NONUNIMODULAR TRANSITIVE GRAPHS

BY ÁDÁM TIMÁR[1]

*Indiana University*

We extend some of the fundamental results about percolation on unimodular nonamenable graphs to nonunimodular graphs. We show that they cannot have infinitely many infinite clusters at critical Bernoulli percolation. In the case of heavy clusters, this result has already been established, but it also follows from one of our results. We give a general necessary condition for nonunimodular graphs to have a phase with infinitely many heavy clusters. We present an invariant spanning tree with $p_c = 1$ on some nonunimodular graph. Such trees cannot exist for nonamenable unimodular graphs. We show a new way of constructing nonunimodular graphs that have properties more peculiar than the ones previously known.

**1. Introduction.** By a *percolation*, we mean some random subgraph of a given graph $G$. We usually assume that the automorphism group of $G$ acts transitively on $G$ and assume (without mentioning this) that the percolation has a distribution that is invariant under this group. The edges (vertices) contained in the random subgraph of the percolation are called *open*, while those not contained in it are called *closed*. A component of the percolation is also called a *cluster*.

An important class of percolations is Bernoulli($p$) *site* (*bond*) *percolation*, where in we remove every vertex (edge) from the graph with probability $1 - p$ and independently of each other.

Research in the past decade has led to a better understanding of percolation on transitive graphs, especially transitive nonamenable graphs. For a survey, see [11] or [15]. However, there is a class of nonamenable graphs where the most important technical tool, the so-called Mass Transport Prin-

Received October 2005; revised May 2006.

[1]Supported in part by NSF Grant DMS-02-31224 and Austrian Science Fund (FWF) P15577.

ciple, fails in its most effective form. Almost none of the general results was possible to extend to this class. These are the *nonunimodular* graphs.

Our knowledge of group-invariant percolation for these graphs has suffered from serious gaps and, as the example in Section 6 shows, surprising phenomena can occur in comparison with unimodular graphs.

One of the fundamental questions about Bernoulli percolation concerns the existence of different phases and behavior at critical points. The critical probability $p_c := \sup\{p : \text{there are no infinite clusters at } p\text{-percolation}\}$ is well known to be equal to $p_u := \inf\{p : \text{there is a single infinite cluster at } p\text{-percolation}\}$ when the underlying graph is amenable, by the arguments in [3] and [5]. On the other hand, $p_c < p_u$ holds for *some* Cayley graph of each nonamenable group [13]. For nonunimodular graphs, a new value of criticality, $p_h$, appears in [7]. To define this, consider a left Haar measure (which is unique up to a constant factor) on the automorphism group of $G$. For each vertex, assign a *weight* that is the Haar measure of its stabilizer. We call a cluster *heavy* if the sum of weights over its vertices is infinite. Otherwise, we call the cluster *light*. (An equivalent, combinatorial, definition is given in Section 2.) Let us mention that a cluster in some Bernoulli percolation is heavy if and only if it contains infinitely many vertices of each weight that is present in $G$, as shown in Corollary 5.10.

Now, let $p_h := \inf\{p : \text{there is some heavy cluster at } p\text{-percolation}\}$. It is easy to show that a unique component is necessarily heavy. The more subtle phases that arise here merit interest, as well as a possible characterization of nonunimodular graphs according to whether or not a phase is present in their case. Finding an example of a graph with $p_c < p_h < p_u < 1$ has been an open problem since [7]. Our result in Section 5 shows that if there are infinitely many heavy clusters, then the intersection of any of them with some finite union of levels contains some infinite *connected* graph. This result makes it simple to rule out the phase of infinitely many heavy clusters in the case of many graphs (a necessary condition for having such a phase is contained in Corollary 5.8), and it may help in finding an example with all possible phases. In Section 3, (3.4), we present a new family of nonunimodular graphs. Most of the previously known examples cannot possibly have $p_c < p_h < p_u < 1$, either because our necessary condition about having infinitely many heavy clusters fails for them or because they have infinitely many ends and hence $p_u = 1$. Our example might be a candidate for a graph with all possible phases.

The question of whether Bernoulli($p_c$) percolation has infinite clusters is one of the central questions in classical percolation theory. A negative answer was given for $\mathbb{Z}^2$ by Harris [9] and Kesten [10] and for $\mathbb{Z}^d$ with $d \geq 19$ by Hara and Slade [8]. For the remaining dimensions, the question has been open for almost fifty years. While the problem for these amenable graphs has thus for resisted all attempts to solve it, for nonamenable unimodular graphs it is known that there are no infinite clusters at criticality; see [1] and [2]. For nonunimodular graphs, this problem has remained unresolved. Lyons, Peres and Schramm have shown (in a paper currently in preparation) that there cannot be infinitely many heavy clusters at criticality. In Section 5 we present this as a corollary of the same result in the unimodular case and our aforementioned theorem concerning infinitely many heavy clusters. The latter uses a key idea from their proof. In Section 4, we show that there cannot be any light clusters at criticality, answering a particular question from [1]. Unfortunately, we cannot rule out the possibility of a unique infinite cluster at criticality, so the problem of critical percolation is not completely solved yet for nonunimodular graphs. For the classes of decorated trees and Diestel–Leader graphs (defined in Section 3) [14] gives a complete solution.

Section 2 gives the definitions of unimodularity, levels, and so on, and introduces the Mass Transport Principle (MTP) and some notation. Section 3 presents the known examples of nonunimodular graphs, together with a new construction that yields graphs with surprising properties. In Section 4, we show that there are no light clusters at critical Bernoulli percolation. In Section 5, we prove that in the case of infinitely many heavy clusters, some subgraph induced by finitely many levels contains some infinite connected open subgraph. One of the corollaries is the nonexistence of infinitely many heavy components at critical percolation (first proved by Lyons, Peres and Schramm). Finally, in Section 6, we present an invariant random spanning tree $T$ on a nonunimodular graph such that $p_c(T) = 1$ a.s. Such an example cannot exist among unimodular graphs; our example answers a question in [1].

**2. Terminology and notation; the MTP.** The graphs we consider are always locally finite.

For a connected graph $G$, fix a vertex $o$ and define $p_c(G) := \inf\{p \in [0, 1]:$ the open component of $o$ at Bernoulli$(p)$ edge percolation is infinite with positive probability$\}$. The *critical probability* $p_c(G)$ is independent of the choice of $o$.

For a graph $G$, we shall always denote its vertex set by $V(G)$, its edge set by $E(G)$ and its group of automorphisms by $\text{Aut}(G)$. The ball of radius $r$ around a vertex $x$ is denoted by $B(x, r)$. Given vertices $x, y \in V(G)$, let $S_x := \{g \in \text{Aut}(G): g(x) = x\}$ be the stabilizer of $x$. Then $S_x y$ is the orbit of $y$ under $S_x$.

A locally compact group is called *unimodular* if its left-invariant Haar measures are also right-invariant. We say that a transitive graph $G$ is unimodular if its group of automorphisms is unimodular. An equivalent definition is the following: $G$ is unimodular if and only if for any $x, y \in V(G)$, we have $|S_x y| = |S_y x|$; see [18]. It is known that amenable graphs are unimodular [16].

A *level* of a nonunimodular graph is a maximal set $X$ of vertices such that for any $x, y \in X$, $|S_x y| = |S_y x|$. Let $G$ be a nonunimodular graph. The *weight* $w(x)$ of a vertex $x$ or of the level of $x$ will be the measure of $S(x)$ by some fixed left-invariant Haar measure. This is equivalent to the following, combinatorial, definition for the weights (which is also determined only modulo a constant factor). Fix some vertex $o$ and define $w(o) := 1$. For each $x \in V(G)$, let $w(x) := |S_x o|/|S_o x|$. This definition is independent of the choice of $o$ up to a constant factor. It is easy to show that for any $x, y \in G$, $|S_x y|/|S_y x| = w(x)/w(y)$. Hence, two vertices are on the same level if and only if their weights are the same. For more details about weights, see [12]. By a *union of levels* we always mean (in a slightly sloppy usage) the subgraph of $G$ induced by the vertices in the levels.

At one point in the paper, we shall refer to *quasi-transitive* graphs: a graph is quasi-transitive if its group of automorphisms has finitely many orbits. Note that if we take a finite set of levels in a transitive nonunimodular group $G$ and restrict $\mathrm{Aut}(G)$ to the union of these levels [i.e., take the subgroup of $\mathrm{Aut}(G)$ of elements that fix levels], it is quasi-transitive and unimodular.

We shall say that a vertex $x$ is *above* vertex $y$ if $|S_x y| > |S_y x|$. If $x$ is above $y$, then $y$ is said to be *below* $x$. We can extend this notion of being "above" to sets: we say that $A$ is above $B$ if every element of $A$ is above any element in $B$, similarly for "below." The *highest* element of a set is one that is above every other element.

We do not usually state explicitly "with probability 1" in the paper when this is the case. By a slight abuse of terminology, but without any real ambiguity, a (*simple*) *path* sometimes means a graph and sometimes a directed graph started from a specified vertex. The length of a path is the number of edges in it.

Given $A$, an event, $A^c$ will stand for its complement. *Insertion tolerance* of a percolation process means that whenever $A$ is some measurable event with $\mathbf{P}[A] > 0$, then $\mathbf{P}[\{\kappa \cup \{e\} : \kappa \in A\}] > 0$. The event $\{\kappa \cup \{e\} : \kappa \in A\}$ was obtained from $A$ by *inserting* $e$. We define *deletion* and *deletion tolerance* analogously, replacing insertion of edges by deletion of them.

The Mass Transport Principle (MTP) is defined for transitive unimodular graphs as follows:

LEMMA 2.1. *Let $G$ be a transitive unimodular graph and $f(x,y)$ a nonnegative function from $V(G) \times V(G)$ that is diagonally invariant under $\mathrm{Aut}(G)$, that is, $f(x,y) = f(\gamma x, \gamma y)$ for any $\gamma \in \mathrm{Aut}(G)$. Then*

$$\sum_{y \in V(G)} f(x, y) = \sum_{y \in V(G)} f(y, x)$$

*for every vertex $x$.*

For a proof, see, for example, [1]. We shall typically define a function $\phi$ on $V(G) \times V(G) \times 2^G$ with $\phi(x, y, \kappa) = \phi(\gamma x, \gamma y, \gamma \kappa)$ for every $\gamma \in \mathrm{Aut}(G)$ and such that $\phi$ is measurable in $\kappa$ with respect to the percolation that we are considering. We refer to this by saying "$x$ sends mass $\phi(x, y, \kappa)$ to $y$ when $\kappa$" (and similarly, "$y$ receives from $x \ldots$"). Then $f(x, y)$ is defined as $\mathbf{E}[\phi(x, y, \kappa)]$. The MTP effectively says that the expected total mass that a vertex sends out is the same as the expected total mass that it receives.

This intuitive statement is no longer true if $G$ is a nonunimodular transitive graph. A simple illustration is the following. Let $G$ be the binary grandmother graph, as in (3.2), and let every vertex send mass 1 to its children. The mass sent out is 2, while the mass received is 1.

The nonunimodular version of the MTP is the following:

LEMMA 2.2. *Let $G$ be a transitive nonunimodular graph and $f(x, y)$ a nonnegative function from $V(G) \times V(G)$ that is diagonally invariant under $\mathrm{Aut}(G)$. Then*

$$\sum_{y \in V(G)} f(x, y) = \frac{1}{w(x)} \sum_{y \in V(G)} f(y, x) w(y)$$

*for every vertex $x$.*

See [1] for more details, proofs and the history of the method.

The "weakness" of the nonunimodular version of the MTP is shown by the following example. Whereas, given a percolation on a unimodular graph, one cannot assign a finite subset of vertices to certain infinite clusters using an invariantly defined property, this is not the case for nonunimodular graphs. Infinite Bernoulli($p$) clusters of a grandmother graph have finitely many uppermost vertices if $p < 1$. The lack of this kind of homogeneity in the vertices of a cluster is responsible for the fact that most known proofs for percolation on unimodular graphs do not extend to nonunimodular ones.

## 3. Examples of nonunimodular graphs.

3.1. *The "Grandmother graph".* Consider a $(k+1)$-regular "family tree," that is, a tree wherein there is some distinguished infinite ray that defines a unique parent and $k$ children for each vertex. Connect every vertex to its grandparent. The resulting graph is called the *grandmother graph*; its levels coincide with generations. This example is due to Trofimov [18]. One can obtain similar examples of nonunimodular graphs by "decorating" a regular tree in some other way to produce a graph whose automorphisms fix a ray of the tree.

3.2. *Diestel–Leader graphs.* The Diestel–Leader graph $DL(k,n)$, first defined in [4], is constructed as follows. Let $T$ and $T'$ be a $(k+1)$-regular and an $(n+1)$-regular tree, respectively, and suppose that their vertices are arranged into "layers" corresponding to the integers, similarly to the tree of the grandmother graph. This is done in such a way that a vertex of $T$ on the $i$th layer has $k$ children on the $(i+1)$th layer and its parent on the $(i-1)$th layer; a vertex of $T'$ on the $i$th layer has $n$ children on the $(i-1)$th layer and the parent on the $(i+1)$th layer. Let the layer of $v$ in $T$ (resp. $T'$) be denoted by $l_T(v)$ [resp. $l_{T'}(v)$]. The vertex set of $DL(k,n)$ will be $\{(x, x') \in V(T) \times V(T') : l_T(x) = l_{T'}(x')\}$. There is an edge between $(x, x')$ and $(y, y')$ iff $x$ and $y$ are adjacent in $T$ and $x'$ and $y'$ are adjacent in $T'$. Hence, a walk on the Diestel–Leader graph can be thought of as two simultaneous walks on the two trees such that the two walkers are always on the same layer. Now, if $k \neq n$, then $DL(k, n)$ is nonunimodular. A level of this nonunimodular graph coincides with the set of vertices whose components are coming from the same layers of the two trees.

Let us recall that $DL(k, k)$ is a unimodular graph which is a Cayley graph of the lamplighter group, with "lamps" being copies of $\mathbb{Z}_k$. For more details, see [19].

3.3. *Free and direct products.* One can obtain new nonunimodular graphs by taking some product of an arbitrary nonunimodular and a unimodular graph. Slightly differing types of direct product can be found in standard graph theory books. The *free product* of two transitive graphs $G_1$ and $G_2$ can be defined inductively as follows. Take a copy of $G_1$ and infinitely many copies of $G_2$. Fix some bijection between the vertices in $G_1$ and the copies of $G_2$ and identify each vertex of $G_1$ with an arbitrary vertex in its image by the bijection. Call the resulting graph $H_1$. Now take a bijection between the vertices $I_1$ of $H_1$ that were not born by identification in the previous step and fix some bijection between $I_1$ and a set of countably many copies of $G_1$. Identify every vertex of $I_1$ with an arbitrary vertex in its image by the bijection to obtain $H_2$. Similarly, given $H_i$, and if $I_i$ is the set of vertices in $H_i$ that were not born by identification in some previous step, take a bijection between $I_i$ and a set of infinitely many copies of $G_j$, where $j \in \{1, 2\}$ is congruent to $i$ mod 2. Identify every vertex of $I_i$ with an arbitrary vertex in its image by the bijection to obtain $H_{i+1}$. Finally, the free product of $G_1$ and $G_2$ is defined as the resulting graph when repeating the above procedure in $\omega_0$ steps. It is easy to see that we get a transitive graph.

Suppose that we take the free product of two nonunimodular graphs where the weights in the two are not the powers of the same constant (e.g., two grandmother graphs of different degrees). We then get an example of a nonunimodular graph where the level structure is not discrete. In other words, the logarithms of the weights in the product graph form a dense subset of $\mathbb{R}$.

3.4. *Constructions from "self-similar" graphs.* The next example (or family of examples) is presented here for the first time. The nonunimodular graph that arises has the property that every level induces a connected graph; there is no previously known example with this feature.

Construct $G$ as follows. Consider $DL(2,2)$, a unimodular graph, as mentioned above. The levels of $G$ will all have weights which are powers of 4.

Each level of $G$ will induce a copy of $DL(2,2)$. One can show that any automorphism of $DL(2,2)$ preserves the layers coming from its construction. (This fact is not essential in what follows, since we could restrict ourselves to only these automorphisms, whose group is obviously transitive.)

Consider two copies of $DL(2,2)$, calling them $G_1$ and $G_2$. We shall construct edges between $G_1$ and $G_2$ and the resulting graph gives the pattern of edges between any two consecutive levels of $G$.

Partition the vertices of $G_2$ into "bags" of four vertices as follows. Let $T$ and $T'$ be the trees that defined $DL(2,2)$, as in (3.2). If $x_1$ and $x_2$ are siblings in $T$ and $y_1$ and $y_2$ are siblings in $T'$, then let a bag be $\{(x_i, y_j) : i, j \in \{1, 2\}\}$. Note that the bags are classes of imprimitivity for the automorphism group of $G_2$, by which we mean that any automorphism preserves these classes. Define a graph $G'_2$ (called the *bag graph*) on the bags: let two bags be adjacent if some element in one of them is adjacent to an element in the other. The graph $G'_2$ is isomorphic to $DL(2,2)$. Furthermore, it has the following properties:

(i) any automorphism of $G_2$ generates an automorphism of $G'_2$,
(ii) conversely, any automorphism of $G'_2$ arises from some automorphism of $G_2$.

Fix some isomorphism $\iota$ between $G'_2$ and $G_1$, and connect a vertex $x \in G_1$ to a vertex $y \in G_2$ if the bag $\iota(x)$ contains $y$.

Now, assign a copy of $DL(2,2)$ to each integer to obtain the graphs $G_i$ ($i \in \mathbb{Z}$) that will be on the levels of our $G$. Define edges between $G_i$ and $G_{i+1}$ as we did between $G_1$ and $G_2$ above. The corresponding isomorphisms can be chosen arbitrarily.

We need to show that the resulting graph is transitive and nonunimodular. For transitivity, consider some automorphism $\gamma$ acting on the graph $G_i$ on some level. By (i), $\gamma$ induces an automorphism on the bag graph $G'_i$ and hence induces an automorphism on the graph $G_{i-1}$. The simultaneous automorphism respects the edges between $G_i$ and $G_{i-1}$. On the other hand, by (ii), any automorphism of the bag graph on $G_{i+1}$ is induced by some automorphism of $G_{i+1}$. Hence, the automorphism that $\gamma$ naturally generates on the bags on $G_{i+1}$ [map the bag whose vertices are adjacent to $x \in G_i$ into the bag whose vertices are adjacent to $y \in G_i$ iff $y = \gamma(x)$] can be extended to an automorphism of $G_{i+1}$. This simultaneous automorphism again respects the edges between $G_i$ and $G_{i+1}$.

We conclude that an arbitrary automorphism $\gamma$ of $G_i$ can be extended to the levels below it and above it and hence recursively to the entire $G$. The same observation is true when $\gamma$ is not an automorphism, but rather an isomorphism from $G_i$ to some $G_j$. Since $\text{Aut}(G_i)$ is transitive, we have that $G$ is transitive.

The verification that $G$ is nonunimodular is effectively just a simple calculation of orbit sizes. If $x \in G_i$ and $y \in G_{i+1}$ are adjacent, then $|S_x y| = 4$ and $|S_y x| = 1$, as can be confirmed by looking at the graph induced by the neighborhood of a vertex in $G$.

One can obtain further examples of nonunimodular graphs by putting some other graph $H$ on the levels instead of $DL(2, 2)$. The graph needs to be such that its vertices can be partitioned into finite classes of imprimitivity ("bags") so that the bag graph $H'$ is isomorphic to $H$ and so it satisfies (i) and (ii) (with $G_2$ and $G'_2$ replaced by $H$ and $H'$, resp.). These are the properties that we referred to as being "self-similar". We have constructed some other such connected graphs, but all those constructions use regular trees in some way. It would be interesting to have other examples of "self-similar" graphs, constructed in different ways.

**4. No light clusters at criticality.** In this section, we prove that Bernoulli($p_c$) edge percolation on a nonunimodular transitive graph does not have infinitely many light infinite clusters. This answers a question from [1]. Lyons, Peres and Schramm have proven (in a paper under preparation) that there cannot be infinitely many heavy clusters at critical percolation. For grandmother graphs and Diestel–Leader graphs, [14] proves that there are only finite clusters at criticality. Unfortunately, the general question as to whether there can be any infinite cluster at Bernoulli($p_c$) edge percolation is still open for nonunimodular graphs.

Given level $\ell$, and level $\ell'$ below $\ell$, the set of vertices above $\ell'$ and below $\ell$ together with $\ell$ (but not including $\ell'$) will be called the vertices *between* $\ell$ and $\ell'$. Let $G(\ell', \ell]$ denote the subgraph of $G$ induced by vertices between $\ell$ and $\ell'$.

We will need a lemma concerning a branching process with bounded interactions. It may be well known, but we have not been able to find a reference. We consider a random-vertex-labeled tree in which the number and labels of the children of a vertex in the $m$th generation $O_m$, given the entire generation $O_m$ (together with the labels there), do not depend on the earlier generations. In the lemma, we will not say explicitly that the vertices are labeled because the label can be considered as encoded in the vertex. We mention the labels here only to emphasize that the distribution of the offspring does not only depend on the size of the current generation, but can also depend on what particular vertices are in it.

Lemma 4.1. *Consider a random rooted tree with the following properties. Fix some $p > 0$ and define $O_0 := \{o\}$, where o is the root. If a generation $O_m$ is already given, then the number of children that the vertices in $O_m$ will have depends only on $O_m$ (and not on the past). Each vertex of $O_m$ has at least $k$ children with probability $\geq p$ and 0 children otherwise. Further, the number of children of a vertex $x \in O_m$ is independent of the number of children of all but at most $\alpha$ of the other vertices in $O_m$. Denote by $O_{m+1}$ the set of children of the vertices in $O_m$.*

*Then the tree is infinite with positive probability whenever $kp > 1$.*

Proof. We may assume that every vertex has exactly $k$ children with probability $p$ and 0 children otherwise. Given a generation $O_m$ with $n$ vertices, let $X_i$ be the random variable defined to be $k$ if the $i$th vertex in $O_m$ has $k$ children and 0 otherwise $(i = 1, \ldots, n)$.

For the second moment of $X := X_1 + \cdots + X_n$, we have

$$\mathbf{E}[X^2] \leq \sum_{i,j=1}^{n} \mathbf{E}[X_i]\mathbf{E}[X_j] + \alpha k^2 n,$$

where the first term here is an upper bound for the summands where $X_i$ and $X_j$ are independent, while the second term is an upper bound for the summands where $X_i$ and $X_j$ are not independent. Define the martingale $Y_m := \frac{|O_m|}{(pk)^m}$. Using the above upper bound and $\mathbf{E}[X_i] = kp$, we have

$$\mathbf{E}[Y_m^2 | Y_{m-1}] \leq (pk)^{-2m}((pk)^2|O_{m-1}|^2 + \alpha k^2 |O_{m-1}|)$$
$$= Y_{m-1}^2 + \alpha k^2 Y_{m-1}(pk)^{-m-1}.$$

Hence,

$$\mathbf{E}[Y_m^2] \leq \mathbf{E}[Y_{m-1}^2] + \alpha k^2 \mathbf{E}[Y_{m-1}](pk)^{-m-1} \leq \alpha k^2 \sum_{i=1}^{m+1} (pk)^{-i}$$

by induction and using $\mathbf{E}[Y_{m-1}] = 1$.

We conclude that there is a uniform upper bound for the second moments of the $Y_m$ because $(kp)^{-1} < 1$. This implies, by the martingale convergence theorem, that the $Y_m$'s converge with probability 1 to a random variable whose mean is the same as that of the $Y_m$'s, namely 1. Hence, this limit variable is positive with positive probability and the $Y_m$'s are all positive on some set of positive probability. On this event, the tree does not die out. □

A *separating set of levels* will be the set $L$ of levels between two fixed levels $l_1$ and $l_2$, $l_2$ below $l_1$, such that there is no path from the vertices above $l_1$ to the vertices below $l_2$ that is disjoint from $L$. In particular, for those nonunimodular graphs where the weight of any level is a power of some constant (so that the graph has a discrete "level structure"), a separating set of levels is just a finite set of levels, coming one below the other. But, in other cases, there are infinitely many levels in $L$, the weights of which form a dense set in some interval.

LEMMA 4.2. *Let $G$ be a nonunimodular transitive graph and $o$ a vertex of $G$. Let $L_0, L_1, \ldots$ be pairwise disjoint separating sets of levels such that $o \in L_0$, $L_j$ is below $L_i$ if $j > i$, and $\bigcup_i L_i$ contains all the vertices below $o$. Consider Bernoulli($p$) edge percolation on $G$ that has light open infinite clusters. Then given the event that $C(o)$ is an infinite light cluster, $|C(o) \cap L_i| \to \infty$ as $i \to \infty$.*

PROOF. By lightness, $|C(o) \cap L_i| \geq 1$ for every $i$ because, otherwise, the cluster would be bounded from below. (We note that a much stronger claim is true for heavy clusters, by a simple argument in the next section.)

Fix $k \in \mathbb{N}$ and let $E(k) = E$ be the event that $C(o)$ is infinite and light and infinitely many of the $L_i$ intersect $C(o)$ in fewer than $k$ elements.

We shall define a mapping for each $j$ from the event $E$. Given a configuration in $E$ such that $|C(o) \cap L_j| < k$, let $C_j$ be the open component of $o$ in the subgraph of $G$ induced by $L_j$ and the levels above $L_j$. Close every open edge that connects a vertex $v$ in $C_j$ to a vertex below $L_j$ (i.e., a vertex in $L_{j+1}$). The probability of the resulting event $E_j$ is $\mathbf{P}[E_j] > c\mathbf{P}[E; |C(o) \cap L_j| < k]$ with some constant $c > 0$ independent of $j$, due to the following argument. A configuration $\kappa$ in $E_j$ could arise from at most a constant number of configurations in $E$: knowing the configuration in $E_j$, we recover the one in $E$ by inserting some of the edges that are adjacent to some vertex in the intersection of $L_j$ and the cluster of $o$ in $\kappa$ (there are $\leq k\delta$ such edges, where $\delta$ is the degree in $G$). On the other hand, $E_j \cap E_i = \varnothing$ for any $i \neq j$ because in the configurations of $E_j$, the cluster of $o$ intersects $L_j$ but does not intersect any $L_i$ for $i > j$. Hence,

$$1 \geq \sum_j \mathbf{P}[E_j] > c \sum_j \mathbf{P}[E; |C(o) \cap L_j| < k]$$
$$= c\mathbf{P}[E]\mathbf{E}[\#L'_j s \text{ such that } |C(o) \cap L_j| < k | E].$$

By the definition of $E$, the right-hand side of the inequality would be infinite if $\mathbf{P}[E] > 0$. So $P[E] = 0$ and the claimed assertion follows. □

THEOREM 4.3. *Let $G$ be a transitive nonunimodular graph. Then there are no infinite light clusters in critical Bernoulli edge percolation on $G$.*

PROOF. Suppose that there exist some light clusters at criticality. Fix a vertex $o$ on level $\ell_0$ and let $\ell_1$ be the lowest level that contains a vertex adjacent (in $G$) to $o$. Let the edges between $o$ and $\ell_1$ be called *long edges* (*from $o$*). The property defined in the next paragraph for a cluster $C(o)$ will allow us to simplify some of our arguments. Fix some $o'$ below $o$ which is connected to $o$ by some long edge. Consider the subgraph of $G$ induced by the vertices on level $\ell_1$ and below it and let $G'(o) = G'$ be the union of this subgraph and the long edge from $o$ to $o'$. One can define $G'(x)$, for any $x \in V(G)$, as the image of $G'(o)$ by some (arbitrarily fixed) automorphism of $G$ that takes $o$ to $x$. Of course, $G'(o)$ depends on the $o'$ we choose and $G'(x)$ depends on both $o'$ and the automorphism chosen.

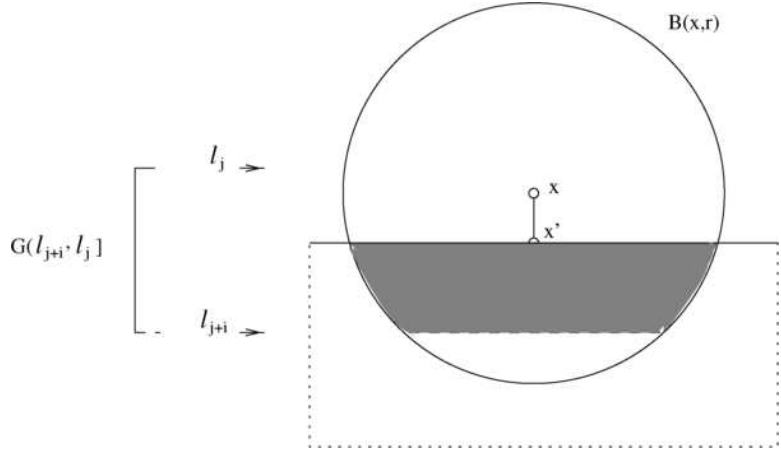

Fig. 1. *Definition of $B_x(i;r)$.*

We say that $C(x)$ is *nice* if $C(x) = C(x)|_{G'(x)}$. Let $F(x)$ be the event that $C(x)$ is infinite, light and nice. Note that the probability of $F(x) =: F$ is independent of $x$ because these events can be mapped into each other by some automorphisms and our probability measure is invariant under automorphisms. By insertion and deletion tolerance and the assumption that there are light clusters, we have $q := \mathbf{P}[F] > 0$. To see this, consider the event that $o'$ is an uppermost vertex of a light cluster, insert the edge between $o$ and $o'$ and delete all other edges incident to $o$. In the resulting configurations (a set of positive probability), $C(o)$ is nice.

We may assume that the weight of $\ell_0$ is 1. If the weight of $\ell_1$ is $\mu$, then let $\ell_i$ be the level with weight $\mu^i$. Note that $(G(\ell_{i+1}, \ell_i])_i$ is a sequence of separating sets of levels that satisfies the properties of the $(L_i)_i$ in Lemma 4.2.

Given an $x$ in some $G(\ell_{j+1}, \ell_j]$, the vertex $x'$ [as described in the definition of $G'(x)$] and $r, i \geq 1$, define $B_x(i; r) := B(x,r) \cap G'(x) \cap G(\ell_{j+i}, \ell_j]$ (see Figure 1). We say that a vertex $v$ in $B_x(i;r) \setminus \{x, x'\}$ belongs to the *side boundary* of $B_x(i; r)$ if there is some edge $\{v, w\}$ such that $w$ is above $\ell_{j+i}$ but not in $B_x(i; r)$. Given $k \in \mathbb{N}$, we say that the open component of $x$ in $B_x(i;r)$ is *good* if it contains at least $k$ vertices in $G(\ell_{j+i}, \ell_{j+i-1}]$ and the side boundary of $B_x(i;r)$ is disjoint from $C$.

From Lemma 4.2, we know that given the event that $C(o)$ is nice and light, with probability arbitrarily close to 1, $|C(o) \cap G(\ell_{i+1}, \ell_i]| \geq k$ if $i$ is sufficiently large. Hence, for any $k$, there exist an $i$ and an $r \in \mathbb{N}$ such that the open component of $o$ in $B_o(i;r)$ is good with probability $> q/2$. (First, we choose $i$, using Lemma 4.2, to satisfy the first property for being good. Then we choose $r$ to satisfy the second property. This can always be done since the intersection of a light cluster with $G(\ell_{j+i}, \ell_j]$ is finite by definition and thus the probability that it intersects a large sphere around $o$ can be made arbitrarily small.) Clearly, there is also a uniform choice, so that $B_x(i;r)$ is good for any $x$ below $o$.

Let $\delta$ be the degree in $G$. Fix $k$ so that $\frac{kq}{2\delta} \geq 4$. Also, fix $i$, $r$ for this $k$ as above [so that the component of $o$ in $B_o(i;r)$ is good with probability $> q/2$]. In what follows, write $B_x$ for the resulting $B_x(i;r)$ $[x \in V(G)]$.

Now, take a $\tilde{p} < p_c$. If $\tilde{p}$ is sufficiently close to $p_c$, then the Bernoulli($\tilde{p}$) component of $x$ in $B_x(i;r)$ is still good with probability $> q/2$, whatever $x$ is. This is true because there is a uniform bound on the sizes of the $B_x$. Fix such a $\tilde{p}$.

Hereafter, we consider Bernoulli($\tilde{p}$) edge percolation; this is where open components, and so on, are understood.

We shall define a branching process $T$ as follows. The vertex $o \in V(G)$ corresponds to the root $\hat{o}$ of $T$. The 0th generation $O_0$ of $T$ is hence $\{\hat{o}\}$. If the open component of $o$ in $B_o$ is good, then it contains at least $k$ vertices in $G(\ell_{i+1}, \ell_i]$. For each of these vertices $x$, add a child $\hat{x}$ to $\hat{o}$ in $T$ and let these constitute the first generation $O_1$ in $T$. If the component of $o$ in $B_o$ is not good, let $\hat{o}$ have 0 children.

Similarly, suppose that we defined the $g$th generation $O_g$ of $T$ using vertices in $G(\ell_{gi+1}, \ell_{gi}]$ so that to each vertex $\hat{x}$ of $O_g$, there corresponds a vertex $x$ in $G(\ell_{gi+1}, \ell_{gi}]$. We can partition $O_g$ so that $\hat{x}$ and $\hat{y}$ are in the same class of the partition iff for the corresponding $x, y \in G(\ell_{gi+1}, \ell_{gi}]$, we have $x' = y'$ (with $x', y'$ defined as in the definition of $G'$). Each set of the partition has $\leq \delta$ elements. Now, choose one vertex in each class uniformly and independently; call the set of chosen vertices *parental* vertices. If $\hat{x} \in O_g$ is not parental, then let it have 0 children.

If $\hat{x}$ is parental, assign $f$ children to it iff the open component of $x$ in $B_x$ is good and contains $f$ vertices from $G(\ell_{(g+1)i+1}, \ell_{(g+1)i}]$; assign 0 children to $\hat{x}$ otherwise. Doing this for each element of $O_g$ defines $O_{g+1}$. Note that $f \geq k$ by being good and that a vertex has $\geq k$ children with probability $> \frac{q}{2}\frac{1}{\delta}$.

The important consequence of the second property for being good is that even though $B_x$ and $B_y$ may intersect each other for certain $x, y \in G(\ell_{gi+1}, \ell_{gi}]$, if the open components of $x$ in $B_x$ and $y$ in $B_y$ are good, then these components are disjoint unless $x' = y'$. This is the case because if $B_x \cap B_y \neq \varnothing$, then there exists an open path between $x$ and $y$ within $B_x \cup B_y$, which can only consist of the edges that connect $x$ and $y$ to $x' = y'$, by definition of being good. On the other hand, the case $x' = y'$ cannot occur if both $x$ and $y$ are parental, by definition. If $\hat{x} \in O_g$ and $\hat{y} \in O_g$ do have children, then they are parental and so the open components of $x$ and $y$ in $B_x$ and $B_y$ are disjoint. Hence, distinct vertices in $T$ correspond to distinct vertices in $C(o)$.

Note that by definition, the edge sets in $B_x$ and $B_y$ are disjoint whenever $x \in G(\ell_{gi+1}, \ell_{gi}], y \in G(\ell_{g'i+1}, \ell_{g'i}]$ and $g \neq g'$. Thus, $\mathbf{P}[|O_{m+1}| = i | O_m, \ldots, O_0] = \mathbf{P}[|O_{m+1}| = i | O_m]$ for any $i$.

Moreover, note that the numbers of children for $\hat{x}$ and $\hat{y}$ are independent if $B_x$ and $B_y$ are disjoint. The number of different $B_y$'s that intersect a fixed $B_x$ is at most $|B(x, 2r)| =: \alpha$, where $r$ is as fixed above.

Hence, Lemma 4.1 applies to the random tree $T$ defined above, with $p = \frac{q}{2\delta}$. The nonextinction of $T$ means that the corresponding vertices in $G$ all belong to one (infinite) open cluster. We conclude that Bernoulli($\tilde{p}$) edge percolation produces infinite clusters. This contradicts the choice $\tilde{p} < p_c$. Hence, there is no infinite light cluster at $p_c$. □

**5. Infinitely many heavy clusters.** Let $\mu$ be the maximum of $w(x)/w(y)$, where $x$ and $y$ are adjacent vertices in $G$. Define the (random) 1-*partition* of $G$ as follows. Choose $U \in [0, 1]$ uniformly at random and let $x$ and $y \in V(G)$ be in the same class of the partition iff $\log_\mu w(x)$ and $\log_\mu w(y)$ are in the same interval of the form $[n + U, n + 1 + U)$, $n \in \mathbb{Z}$. Note that the 1-partition of $G$ also partitions its levels. The idea of such a partition comes from Lyons, Peres and Schramm, and it is essential in what follows. We use it to partition the levels similarly to the previous section, when we partitioned the graph to separating sets of levels. The main difference is that here, we shall get a partition that is automorphism-invariant.

The next lemma is similar to one due to Häggström [6].

LEMMA 5.1. *Suppose that $G$ is a nonunimodular transitive graph and consider the 1-partition of $G$ and some other random process $\hat{\omega}$ that is invariant under* $\mathrm{Aut}(G)$ *and independent of the 1-partition. Let $H$ be a random graph on $G$ that is an equivariant function of the 1-partition and $\hat{\omega}$ and such that (*i*) the endpoints of any edge of $H$ are in the same class of the 1-partition and (*ii*) every component of $H$ is a finite tree. Then the expected degree of a vertex, given that it is in $H$, is at most $2\mu$.*

PROOF. Perform the following mass transport. Start with $\deg_H x$ mass in each vertex and redistribute it equally among the vertices in its (finite) cluster. The mass received by a vertex $x$ is the average degree in its finite component. Since the average degree in a finite tree is $< 2$, we get that

$$\mathbf{E}[\deg_H x] = \mathbf{E}[\text{mass sent out}] \leq \mu \mathbf{E}[\text{mass received}] \leq 2\mu,$$

using Lemma 2.2. □

We shall call a set of the form $L = \{x : \mu^a < w(x) \leq \mu^b\}$ a *slab* if $a, b \in \mathbb{R}$ and $b - a \geq 1$. Every separating set of levels (in the previous section) is a slab, but the converse is not true: a slab may not have an uppermost level. It is clear that a slab separates the levels below it from the levels above it. Note that each class in the 1-partition is a slab.

LEMMA 5.2. *A heavy cluster intersects every slab in infinitely many vertices.*

PROOF. Otherwise, by deletion tolerance, with positive probability there would be a heavy cluster that does not have any vertex below (resp. above) a certain level. Then by insertion tolerance, there would also be a heavy cluster with a single lowest (resp. uppermost) vertex. Let every vertex of the cluster send unit mass to this vertex. This contradicts the MTP. □

We call a vertex $x$ in a connected graph $C$ an *encounter point* if $C \setminus x$ has at least three heavy infinite connected components. The next lemma is a version of a similar lemma in [2], with the same proof.

LEMMA 5.3. *Let $\omega$ be the open graph in a Bernoulli edge percolation. Suppose that there are infinitely many heavy clusters in $\omega$. Then every heavy cluster in $\omega$ contains infinitely many encounter points. Furthermore, let $L_0$ be a class of the 1-partition of $G$ and suppose that $L_0 \cap \omega$ contains some encounter point. Then for any encounter point $x \in L_0 \cap \omega$ and any open heavy component $C$ in $\omega \setminus \{x\}$ that is adjacent to $x$ in $\omega$, $C \cap L_0$ contains some vertex that is an encounter point for $\omega$.*

PROOF. The existence of encounter points follows from insertion tolerance. Suppose now that there exists an encounter point $x$ and that $C$ is a heavy component of $\omega \setminus \{x\}$ that is adjacent to $x$ and such that $C \cap L_0$ does not contain any encounter points. Then let each vertex $y$ contained in such a $C \cap L_0$ send mass 1 to $x$ in $L_0 \cap \omega$, where $x$ is an encounter point and where there is an $x$-$y$ path in $\omega$ that does not contain any other encounter point inside $L_0$. Since $C \cap L_0$ is infinite (by deletion tolerance and Lemma 5.2), the expected mass received is infinite. The expected mass sent out is at most 1. This MTP contradiction proves the second part of the claimed assertion. □

We shall present one more proposition before proceeding to the main result of this section. Consider the set

$$R := \{\log_\mu w(\ell) : \ell \text{ is a level of } G\}.$$

Note that $R$ is an additive subgroup of $\mathbb{R}$, generated by the elements of $\{\log_\mu \frac{w(x)}{w(y)} : x \text{ and } y \text{ are adjacent}\}$, which is finite because $G$ is locally finite. Since $R$ is torsion-free, it is isomorphic to $\mathbb{Z}^n$ for some $n$. Fix some isomorphism. Note that any automorphism of $G$ acts on the weights of the levels by multiplying them with a constant, so it acts on $R$ by adding some constant. Since $\text{Aut}(G)$ acts on the set of levels of $G$ transitively, it also acts on $R$ and $\mathbb{Z}^n$ transitively, where this action is "transferred" by the isomorphism between $R$ and $\mathbb{Z}^n$.

PROPOSITION 5.4. *There exists an invariant random exhaustion of the levels of $G$ by a sequence of finite partitions, that is, there exists an invariant random sequence $(\mathcal{P}_i)_i$ such that for every $i$, the $\mathcal{P}_i$ partitions the set of levels of $G$ into finite sets and any two levels are in the same class of $\mathcal{P}_i$ with probability tending to $1$ as $i \to \infty$.*

PROOF. It is enough to prove that there exists an invariant exhaustion for $\mathbb{Z}^n$ by partitions consisting of finite classes (later called *finite partitions*). Then the isomorphism between $R$ and $\mathbb{Z}^n$ will transfer the exhaustion to the levels.

Denote the $i$th coordinate of a vertex $u \in \mathbb{Z}^n$ by $u_i$. Let $\lfloor t \rfloor$ stand for the integer part of a number $t$.

Let $\mathcal{P}_m$ consist of the following classes. Choose an $x \in \{1, 2, \ldots, m\}^n$ uniformly at random and let vertices $v, w \in Z^n$ be in the same class of the partition iff $\lfloor \frac{v_i - x_i}{m} \rfloor = \lfloor \frac{w_i - x_i}{m} \rfloor$ for every $i \in \{1, \ldots, n\}$. It is easy to see that the sequence of partitions defined this way satisfies the required properties. □

THEOREM 5.5. *Consider Bernoulli percolation on some nonunimodular transitive graph $G$ and suppose that it has infinitely many heavy components. Then with probability $1$, for any infinite component $C$, there is some finite union $L$ of levels such that $L \cap C$ has some infinite connected component.*

PROOF. Suppose that for some heavy cluster $C$ and any $L$ that is a finite union of levels, all connected components in $C \cap L$ are finite. Let $\omega$ be the subgraph consisting of these $C$'s. Note that $\omega$ is invariant.

Very briefly, we shall define graphs on the classes of the 1-partition. Hence, in each connected component of this graph, the weights of the vertices will have bounded logarithms. These graphs will be equivariant functions of $\omega$, the 1-partition and some additional randomness. Each of these graphs will be a forest with high enough degrees. The vertices in one tree of the forest are all in the same $\omega$-component. Then we shall construct an invariant exhaustion of the forests by subforests of only finite components, using $\omega$ and Proposition 5.4. This will contradict Lemma 5.1.

Consider the 1-partition on $G$. For each equivalence class $L_0$ of the 1-partition, we shall construct a forest $\Phi(L_0)$ on $L_0$, in a sequence of steps.

By Lemma 5.3, there are encounter points in $\omega$, hence any $L_0$ contains an encounter point with positive probability (since a 1-partition has countably many classes $L_0$ and they are translates of each other by some automorphism). Define a graph $\Gamma_W$ on the set $W$ of encounter points of $\omega$ in $L_0$: put an edge between two of them if there is an open path in $\omega$ between them such that no other point of $W$ is in this open path. Every element of $W$ is an encounter point of $\Gamma_W$ (by which we mean that its deletion results in at least three new infinite components), by Lemma 5.3. Now, define a subforest $F$ of $\Gamma_W$. For every $v \in W$ and each component $I$ of $\Gamma_W \setminus v$ such that $I$ is adjacent to $v$ in $\Gamma_W$, choose uniformly a vertex of $I$ that is closest to $v$ in $\omega$. Put a directed edge from $v$ to this vertex. Doing this for every $v \in W$ and $I$ as above, we obtain a digraph $\vec{M}$. Denote by $M$ the graph that result from ignoring the directions of the edges in $\vec{M}$. There may be cycles in $M$, but any two cycles can share at most a vertex. This is the case for the following reason. Any subgraph $X$ of $\vec{M}$ that is 2-connected in $M$ is obviously contained in a 2-connected component (in other words, a block) of $M$. Hence, the outdegree of a vertex can be at most one in $X$. The union of two cycles sharing more than one vertex would be a 2-connected graph and, of course, it has average degree greater than two. But, then, by the condition on the outdegree, we would also have that the average degree in the graph is at most two, giving a contradiction. So, we can delete a uniformly chosen edge from each of the possibly arising pairwise edge-disjoint cycles in $M$ to obtain a forest $F$. This $F = F(L_0)$ is defined for each $L_0$ with positive probability. Of course, by ergodicity, the constructed forest is almost surely nonempty.

The family of forests $\{F = F(L_0) : L_0 \text{ is a class of the 1-partition}\}$ was constructed in a way that is an equivariant function of the 1-partition, the percolation and some additional randomness that was also invariant (when breaking ties). It has only infinite components and every point $x$ in $F$ has degree $\geq 3$ because each element of $W$ is incident to at least one edge for each infinite component that results from the component of $x$ after the deletion of $x$.

Once we have the $F(L_0)$'s, we can define a forest $\Phi(L_0)$ on the vertices of each class $L_0$ of the 1-partition so that the expected degree of a vertex, given that it is in $\Phi$, is greater than $2\mu$. Choose an invariant partition of $F$ into finite bags so that the expected number of bags that a bag on a fixed vertex is adjacent to is $> 2\mu$. This can be done by means of an "extra" percolation on $V(F)$ with a sufficiently low density of open vertices and then by defining a bag as the set of vertices that are closest to a particular open vertex in the extra percolation (in case of ties, decide by a random uniform choice). Then we can choose the "leader" of this bag to be the open vertex of the bag by the extra percolation. Connect the leader of each bag to the leaders of the adjacent bags. We defined a forest $\Phi(L_0)$ on the leaders and the expected degree given that a vertex is in this forest is $> 2\mu$. Define $\Phi := \bigcup_{L_0} \Phi(L_0)$.

The forest $\Phi$ that we have obtained is an equivariant function of the 1-partition, the percolation and the (independent) extra random variables that we were using. Every vertex has expected degree $> 2\mu$ on that it is in $\Phi$. On the other hand, it has an equivariant exhaustion $\mathcal{R}_i$ by graphs of only finite components, for the following reason. Consider the invariant exhaustion $\mathcal{P}_i$ of the set of levels of $G$ by finite partitions, as in Proposition 5.4. Each class $A$ in $\mathcal{P}_i$ contains finitely many levels of $G$: call their union $L_A$. Consider the connected subgraphs $\omega_A$ of $\omega$ induced by the $L_A$'s as $A$ goes through all the classes in $\mathcal{P}_i$. By our assumption, every $\omega_A$ has only finite components. Now, define $\mathcal{R}_i$ to be the partition that the $\omega_A$'s induce on the vertices of $\Phi$ ($A \in \mathcal{P}_i$): two vertices are in the same class of $\mathcal{R}_i$ iff they belong to the same connected component of some $\omega_A$. Having defined the vertices in a class of $\mathcal{R}_i$ this way, let the graph on them be the one that they induce in $\Phi$.

Each resulting $\mathcal{R}_i$ has only finite connected components and any edge in these graphs has both endpoints in the same class of the 1-partition. So, Lemma 5.1 tells us that the expected degree in the graph $\mathcal{R}_i$ is $< 2\mu$. On the other hand, the sequence $(\mathcal{R}_i)_i$ exhausts $\Phi$ because any set of finitely many levels of $G$ is contained in the same class of $\mathcal{P}_i$ with probability tending to 1 and hence the endpoints of any edge of $\Phi$ are in the same component of $\mathcal{R}_i$ with probability tending to 1. This contradicts the fact that the expected degree in $\Phi$ is $> 2\mu$. □

COROLLARY 5.6. *Suppose that the assumptions of Theorem 5.5 hold. Then there exists some finite union $L$ of levels which induces infinitely many infinite open components.*

PROOF. Apply Theorem 5.5 to two infinite clusters: there are two finite sets of levels, $L_1$ and $L_2$, which both induce some infinite connected open component. Since there are only countably many pairs of finite sets of levels, we can choose $L_1$ and $L_2$ so that, with positive probability, they both induce some infinite connected open component. Let $L_3$ be a finite set of levels such that $L_1 \cup L_2 \cup L_3 := L$ is connected ($L_3$ can be chosen as the set of levels that a finite path that intersects each level of $L_1 \cup L_2$ visits). Then $L$ is connected, $\mathrm{Aut}(G)$ acts quasi-transitively on it and the Bernoulli percolation considered has at least two infinite components, by our assumption on $L_1$ and $L_2$. Then it is well known that by quasi-transitivity, there are infinitely many infinite open components. □

There follow some important corollaries. The first one was proven by Lyons, Peres and Schramm in a paper currently in preparation.

COROLLARY 5.7. *Let $G$ be a transitive nonunimodular graph. Then there cannot exist infinitely many infinite clusters for* Bernoulli$(p_c)$ *percolation on $G$.*

PROOF. We have already established this for light clusters. For heavy clusters, Corollary 5.6 shows that if there were infinitely many of them, then the percolation restricted to some union $L$ of finitely many levels would also have infinitely many infinite clusters. However, the action of Aut$(G)$ on $L$ is quasi-transitive and unimodular, in which case it is known that the existence of infinitely many infinite Bernoulli clusters implies that their critical probability is $< 1$ (see [1]). This would contradict the definition of $p_c$. □

COROLLARY 5.8. *Let $G$ be a transitive nonunimodular graph and suppose that the restriction of $G$ to any finite union of its levels induces only finite components. Or suppose, more generally, that the subgraph induced by any finite set of levels is amenable. Then there cannot be infinitely many heavy clusters for any Bernoulli percolation on $G$.*

PROOF. By Corollary 5.6, this would contradict the assumption on $G$ because an amenable quasi-transitive graph cannot have infinitely many infinite clusters at Bernoulli percolation. □

REMARK 5.9. There do exist nonunimodular graphs with infinitely many heavy clusters. An example is the free product of some arbitrary nonunimodular graph with the $r$-regular tree, where Bernoulli$(p)$ percolation is taken with $p > \frac{1}{r-1}$.

Theorem 5.5 cannot be strengthened to say that one level necessarily contains an infinite open subgraph. Consider the example in the previous paragraph and replace every edge $\{v_0, v_1\}$ within a level by the set of edges that connect $v_i$ to some neighbor of $v_{1-i}$ on some lower level, $i = 0, 1$. For $p$ sufficiently large, this graph will have infinitely many heavy clusters, but any level induces only isolated vertices.

COROLLARY 5.10. *A cluster of some Bernoulli percolation is heavy if and only if it intersects each level in infinitely many vertices a.s.*

PROOF. That infinite intersection implies heaviness is trivial.

If there is a unique cluster, the claim is obvious by the MTP. (Otherwise, let each vertex on a level send mass $1/k$ to the $k$ vertices of the cluster on that level whenever $1 \leq k < \infty$. One may assume $1 \leq k$ by insertion tolerance.) By Theorem 5.5, a heavy cluster $C$ intersects some level $\ell$ in infinitely many vertices. Then it intersects any adjacent level in infinitely many vertices. (Otherwise, similarly to the previous arguments, let each vertex in $C \cap \ell$ send mass $1/|C \cap \ell'|$ to every vertex of $C \cap \ell'$, where $\ell'$ is an adjacent level and where we assume $|C \cap \ell'| \geq 1$. The existence of such an $\ell'$ follows from the assumption and insertion tolerance and we arrive at a MTP contradiction.) Similarly, we obtain the claimed assertion for any level. □

REMARK 5.11. All our arguments work if we only assume insertion and deletion tolerance of the percolation process, instead of Bernoullicity. For the case where we do not have tolerances with some uniform constant, some standard extra argument is needed (see, e.g., [17]). The proofs can be repeated for site percolation instead of bond percolation.

**6. A connected $\omega$ with $p_c(\omega) = 1$ on a nonamenable $G$.** For unimodular graphs, the existence of an invariant connected subgraph of critical probability 1 implies amenability and the existence of an automorphism-invariant spanning tree with at most two ends, as shown in Theorem 5.3 of [1]. We shall define a graph $G$ with infinitely many ends and a nonunimodular group $\mathrm{Aut}(G)$ acting on it transitively, as well as random invariant subgraph $\omega$ of $G$ such that (i) $\omega$ is connected and (ii) for any $p \in (0, 1)$, Bernoulli($p$) percolation on $\omega$ results in only finite components [i.e., $p_c(\omega) = 1$]. This construction answers the question after Theorem 5.3 in [1] negatively because it shows that the existence of such an $\omega$ does not imply the existence of an $\mathrm{Aut}(G)$-invariant spanning tree of $G$ with at most two ends (such a spanning tree cannot exist in a graph with infinitely many ends). In particular, it is an example of a transitive nonamenable graph with an invariant connected (spanning) subgraph of critical probability 1.

Before giving $\omega$, we shall need the following random tree construction. Denote by $\Upsilon_n(x)$ the tree rooted at $x$ that consists of a path $P$ of length $2^n$ starting from the root and then a binary tree of depth $n$ rooted at the other endpoint of $P$. Define a random tree $\Upsilon$ rooted at $o$ as follows. Take Bernoulli($q$) site percolation $\mathcal{B}$ on the vertices of $\mathbb{N}^+$. Let $x_0 = 0$ and $x_1, x_2, \ldots$ be the closed vertices of $\mathbb{N}^+$, in increasing order. We build $\Upsilon$ in consecutive steps. In the 0th step, we have the root $o$. In the $j$th step, we add a tree $\Upsilon_{x_j - x_{j-1} - 1}(z)$ to each leaf $z$ of the tree constructed in step $j - 1$. (When $j = 1$, $o$ is regarded as a leaf.)

The next proposition originally had a lengthier proof. The following short one was suggested by Lyons, Peres and Pete.

PROPOSITION 6.1. *Let $\Upsilon$ be the random rooted tree as above and let $q \geq 1/2$. Then $p_c(\Upsilon) = 1$ almost surely.*

PROOF. Since $\Upsilon$ is spherically symmetric, we have $p_c^{-1}(\Upsilon) = \mathrm{br}(\Upsilon) = \mathrm{gr}(\Upsilon)$ (see, e.g., [12]), where br denotes branching number and gr denotes growth rate. So, it is enough to prove that $\Upsilon$ has subexponential growth. Let $X_i := x_i - x_{i-1} - 1$ be the depth of the $i$th "block". Then the $X_i$'s are i.i.d. with geometric distribution of parameter $1 - q$. After $n$ blocks, a level has size $S_n = 2^{X_1 + \cdots + X_n}$ and it is the $L_n = (2^{X_1} + \cdots + 2^{X_n})$th level in the tree, both understood up to a bounded factor. The log of the lower growth of the tree is then at most

$$\liminf \frac{\log S_{n-1}}{L_n} = \liminf C \frac{(X_1 + \cdots + X_{n-1})}{2^{X_1} + \cdots + 2^{X_n}}.$$

By the law of large numbers, this tends to zero when $q > 1/2$ because in that case, $2^{X_i}$ has infinite mean. □

Denote by $T_k$ the rooted binary tree of depth $k$.

For two graphs $H$ and $H'$, define $H \boxtimes H'$ as the graph on the set $V(H) \times V(H')$, where $(x, x')$ is adjacent to $(y, y') \neq (x, x')$ iff $x$ is adjacent or identical to $y$ in $H$ and $x'$ is adjacent or identical to $y'$ in $H'$. Let $M$ be the grandmother graph constructed from a 3-regular tree. Define $G$ to be the graph $M \boxtimes K_4$, where $K_4$ is the complete graph on four vertices. So, each vertex of $G$ can be written as $(x, j)$, where $x \in V(M)$ and $j \in \{0, 1, 2, 3\}$; for simpler notation, we equip $\{0, 1, 2, 3\}$ with modulo 4 addition. We shall give a random spanning tree $\omega$ on $G$ that is invariant under $\mathrm{Aut}(G)$ and which satisfies $p_c(\omega) = 1$. The tree that we construct will contain edges only from the subgraph $G' := T \boxtimes K_4$ of $G$, where $T$ is the underlying 3-regular tree that we used in the construction in $M$. Of course, $\omega$ will be invariant only under $\mathrm{Aut}(G)$ and not under $\mathrm{Aut}(G')$.

Fix some bijection between the levels of $G$ and the edges of $\mathbb{Z}$. So, we also have a fixed natural bijection between the vertices of $\mathbb{Z}$ and sets of edges in $G'$ that connect two consecutive levels. Consider the following random $\omega$ on $G'$.

First, let $\mathcal{B}$ be the Bernoulli($q$) percolation on $V(\mathbb{Z})$ (with $q$ specified later). If we delete the edges of $G'$ that correspond to closed vertices of $\mathbb{Z}$, then $G'$ breaks apart into finite components; call the set of these components $\mathcal{C}$. Each element of $\mathcal{C}$ has the form $F \boxtimes K_4$, where $F$ is isomorphic to some $T_k$, the root of $T_k$ corresponding to the highest vertex of $F$.

Let $\mathcal{P}_k$ be the set of closed paths on $T_k$ with the properties that (i) they start from the root (and end there) and (ii) that they visit each vertex, but visit each of them at most three times. Since a usual depth-first search results in an element of $\mathcal{P}_k$, it is not empty.

For each component $C \in \mathcal{C}$, first choose an element $j(C) = j$ of $\{0, 1, 2, 3\}$ uniformly at random. If $C \cong T_k \boxtimes K_4$, take an element $Q_0$ of $\mathcal{P}_k$ uniformly at random and define $Q$ to be the simple path on vertices in $C \cap \{(x, i) : x \in V(T), i \in \{j, j+1, j+2\}\}$ that has vertex $(x, i)$ at the $t$th place iff the $t$th vertex of $Q_0$ is $x$ and where this is the $(i - j + 1)$th occurrence of $x$ on $Q_0$. Having defined the subsequent vertices of $Q$ so that adjacent vertices follow each other, the edges of $Q$ are also defined. We shall call the first vertex of $Q$ the *root* of the component $C$. To finish our construction of $C \cap \omega$, define $R$ to be the binary tree $C \cap \{(x, j+3) : x \in V(T)\}$. The last vertex of $Q$ and the highest vertex of $R$ are connected by an edge; we define $C \cap \omega$ to be the tree formed by this edge and by $Q \cup R$.

Finally, if $C'$ is another component in $\mathcal{C}$ that is adjacent to $C$ and below it, then the root of $C'$ is adjacent to four vertices of $C$, say $(v, 0), (v, 1), (v, 2), (v, 3)$. Define the edge between the root of $C'$ and $(v, j(C) + 3)$ to be in $\omega$. Note that $(v, j(C) + 3)$ is a leaf of $C \cap \omega$. Now, from every 4-tuple $\{(x, 0), \ldots, (x, 3)\}$, at least one element belongs to the tree defined thus far—let $(x, i)$ be a uniformly randomly chosen one of them. If some other $(x, i')$ does not yet belong to the tree, then we add the edge between $(x, i)$ and $(x, i')$ to the tree.

The $\omega$ just defined is a spanning tree. Let $X_o$ be the infinite component of $\omega \setminus o$ that intersects the levels above $o$ in an infinite set. Define $\omega_o := \omega \setminus X_o$. For a fixed $o \in V(G)$, consider the event $E_o$ that the edges upward from $o$ are not in $\omega$. On $E_o$, the graph $\omega_o$ can be obtained by subdividing a random tree $\Upsilon$ of Proposition 6.1, with $\mathcal{B}|_{\mathbb{N}^+}$ here playing the role of $\mathcal{B}$ there. Subdivision does not increase the probability that a fixed vertex is in an infinite component after percolation. Hence, on $E_o$, if $p \in (0, 1)$, then the connected component of $o$ after Bernoulli($p$) percolation on $\omega$ (or $\omega_o$) is almost surely finite, by Proposition 6.1. Any cluster in almost every configuration is in $E_o$ for some $o \in V(G)$: otherwise, it would contain some infinite path that has only finitely many vertices below any fixed level. But, there are only countably many such paths in $\omega$ (and each of them contains some closed edge with probability 1).

Hence, with probability 1, every cluster is in some $E_o$, so it contains a vertex that is in a finite cluster. Hence, every cluster is finite. This establishes that $p_c(\omega) = 1$.

**Acknowledgments.** I am indebted to Russell Lyons for helpful discussions and corrections while the drafts were prepared. I also wish to thank Yuval Peres, Gábor Pete and the anonymous referees for their comments on the manuscript.

Department of Mathematics
Indiana University
Bloomington, Indiana 47405-5701
USA
E-mail: [atimar@indiana.edu](atimar@indiana.edu)
URL: [http://mypage.iu.edu/~atimar/](http://mypage.iu.edu/~atimar/)